\documentclass[12pt,twoside]{article}
\usepackage{amssymb,latexsym,amsxtra,amscd,amsfonts,amsthm,amsmath,verbatim}

\makeatletter

\def\section{\@ifstar\unnumberedsection\numberedsection}
\def\numberedsection{\@ifnextchar[
  \numberedsectionwithtwoarguments\numberedsectionwithoneargument}
\def\unnumberedsection{\@ifnextchar[
  \unnumberedsectionwithtwoarguments\unnumberedsectionwithoneargument}
\def\numberedsectionwithoneargument#1{\numberedsectionwithtwoarguments[#1]{#1}}
\def\unnumberedsectionwithoneargument#1{\unnumberedsectionwithtwoarguments[#1]{#1}}
\def\numberedsectionwithtwoarguments[#1]#2{%
  \ifhmode\par\fi
  \removelastskip
  \vskip 5ex\goodbreak
  \refstepcounter{section}%
  \hbox to \hsize{\hss\vbox{\advance\hsize by 1cm
      \hspace*{0.10in}
      \leavevmode\Large\bfseries\raggedright
      \thesection.\
      #2\par
      \vskip 0ex
      }}\nobreak
  \vskip 2ex\nobreak
  \addcontentsline{toc}{section}{%
    \protect\numberline{\thesection}%
    #1}%
  }
\def\unnumberedsectionwithtwoarguments[#1]#2{%
  \ifhmode\par\fi
  \removelastskip
  \vskip 3ex\goodbreak
  \hbox to \hsize{\hss\vbox{\advance\hsize by 1cm
     \hspace*{0.10in}
      \leavevmode\Large\bfseries\raggedright
      #2\par
      \vskip 0ex
            }}\nobreak
  \vskip 0ex\nobreak
  \addcontentsline{toc}{section}{%
    #1}%
  }

\newcommand{\ps@thesis}{
\renewcommand{\@oddhead}{\hfil {\large Symmetric \,Algebra \,For \,Monomial \,Curves}\hfil {\small \thepage}}%
\renewcommand{\@evenhead}{\hfil{\large Debasish \, Mukhopadhyay}\hfil {\small \thepage}}%
\renewcommand{\@oddfoot}{}%
\renewcommand{\@evenfoot}{}%
}

\makeatother

\newtheorem{theorem}{Theorem}[section]

\newtheorem{lemma}[theorem]{Lemma}
\newtheorem{Key-lemma}[theorem]{Key-Lemma}
\newtheorem{corollary}[theorem]{Corollary}

\newtheorem{remark}[theorem]{Remark}

\newtheorem{definition}[theorem]{Definition}

\renewcommand*{\abstract}{}
\renewcommand*{\abstract}{{\noindent \upshape\textbf{Abstract.\,\,}}}

\newcommand*{\Kwords}{}
\renewcommand*{\Kwords}{{\noindent \upshape\textbf{Keywords\,:\quad }}}

\renewcommand*{\proof}{}
\renewcommand*{\proof}{{\noindent \upshape\textbf{Proof.\quad}}}

\newcommand*{\Acknow}{}
\renewcommand*{\Acknow}{{\noindent \upshape\textbf{Acknowledgements.\quad}}}

\begin{document}

\title{\Large\bf{The Symmetric Algebra For Certain Monomial Curves}}
\author{Debasish  Mukhopadhyay}

\date{}

\maketitle

\thispagestyle{empty}

\abstract \begin{small}Let $p \geq 2$ and
$0<m_{0}<m_{1}<\ldots<m_{p}$ be a sequence of positive integers such
that they form a minimal arithmetic sequence. Let $\wp$ denote the
defining ideal of the monomial curve $\mathcal{C}$ in
$\mathbb{A}_{K}^{p+1}$, defined by the parametrization $X_{i} =
T^{m_{i}} \,\, {\rm for} \,\,i \in [0,p] $. Let $R$ denote the
polynomial ring $K[X_{1}, \ldots, X_{p},X_{0}]$. In this article, we
construct a minimal Gr\"obner basis for the symmetric algebra for
such curves, as an $R$-module and what is interesting is that the
proof does not require any $S$-polynomial computation.\end{small}
\medskip

\Kwords {\small Monomial Curves\,,\, Gr\"{o}bner Basis\,,\,
Symmetric
Algebra\,.}\\[2mm]
\section{Notation}
\noindent Let $\mathbb{N}$ denote the set of non-negative integers
and the symbols\, $a,b,d,i,i',j,j',l,l',m,n,p,q,s$ \,denote\\[2mm]
non-negative
integers\,. For our convenience we define \quad $[a\,,\,b]=\{i \,\mid \,a\leq i\leq b\}$ \quad ,\\[2mm]
\noindent $\epsilon(i\,,\,j)=
\begin{cases}
i+j \hspace*{0.25in} {\rm if} \hspace*{0.25in} i+j < p\\
p \hspace*{0.50in} {\rm if} \hspace*{0.25in} i+j \geq p\\
\end{cases}$
\quad and \quad $\tau(i\,,\,j)=
\begin{cases}
0 \hspace*{0.25in} {\rm if} \hspace*{0.25in} i+j < p\\
p \hspace*{0.25in} {\rm if} \hspace*{0.25in} i+j \geq p\\
\end{cases}$
\section{Introduction}
\noindent A class of rings, collectively designated  {\it Blowup
Algebras}, appear in many constructions in Commutative Algebra and
Algebraic Geometry. The ancestor of the blowup algebra is the {\it
symmetric algebra}. Given an $R$-module $M$, the {\it symmetric
algebra} of $M$ is an $R$-algebra $S(M)$ which together with an
$R$-module homomorphism
$$\pi : M \longrightarrow S(M)$$
solves the following universal problem: For a commutative
$R$-algebra $B$ and any $R$-module homomorphism $\varphi : M
\longrightarrow B$, there exists a unique $R$-algebra homomorphism
$\overline{\varphi} : S(M) \rightarrow B$ such that $\varphi =
\overline{\varphi} \circ \pi$. Thus, if $M$ is a free module then
$S(M)$ is the polynomial ring $R[T_{1}, \ldots , T_{m}]$, where $m =
{\rm rank}(M)$. More generally, when $M$ is given by the
presentation
$$ R^{r}\stackrel{\varphi}\longrightarrow
R^{m}\longrightarrow M \longrightarrow 0\,, \quad \varphi =
(a_{ij}),$$ its symmetric algebra is the quotient of the polynomial
ring $R[T_{1}, \ldots , T_{m}]$ by the ideal generated by the
$1$-forms
$$f_{j} = a_{1j}T_{1} + \cdots + a_{mj}T_{m}, \quad j=1,\ldots, r.$$
Conversely, any quotient ring of the polynomial ring $R[T_{1},
\ldots , T_{m}]/\mathcal{L}$, with $\mathcal{L}$ generated by the
$1$-forms in the $T_{i}$'s is the symmetric algebra of a module.
\medskip

\section{Monomial Curves}
\noindent Let $p \geq 2$ and $0<m_{0} < m_{1} < \ldots < m_{p} $ is
an arithmetic sequence of integers with $m_{i}=m_{0}+id$ for
$i\in[1,p]$ and $d \geq 1$. We also assume that $m_{0}=ap+b$ with $a
\geq 1$ and $b \in [1,p]$\,. Let $\Gamma $ denote the numerical
semigroup generated by \,$m_{0},\ldots,m_{p}$\,
i.e.,\,$\Gamma:=\sum_{i=0}^{p} \mathbb{N} m_{i}$. We further assume
that $\gcd(m_{0},d) = 1$ and the set $S=\{m_{0},\ldots,m_{p}\}$
forms a minimal set of generators for $\Gamma$.

\noindent Let $K$ be a field and $X_{1},\ldots, X_{p},X_{0}, T$ are
indeterminates. Let $\wp$ denote the kernel of the $K$-algebra
homomorphism $\eta : R:= K[ X_{1}, \ldots, X_{p}, X_{0} ] \to K[T]$,
defined by $\eta(X_{i}) = T^{m_{i}} $\, for $i\in[0,p]$\,. The prime
ideal $\wp$ is an one-dimensional perfect ideal and it is the
defining ideal of the affine monomial curve given by the
parametrization $X_{i} = T^{m_{i}}\,\,  {\rm for}\,\, i\in[0,p]$. It
is easy to verify that $\wp$ is generated by binomials of the form
\,\,$\displaystyle{X_{1}^{\alpha_{1}}\cdots
X_{p}^{\alpha_{p}}\,X_{0}^{\alpha_{0}}\, - \,X_{1}^{\beta_{1}}\cdots
X_{p}^{\beta_{p}}\,X_{0}^{\beta_{0}}\,\, {\rm with}
\,\,\sum_{i=0}^{p}\alpha_{i}m_{i} =
\sum_{i=0}^{p}\beta_{i}m_{i}}$\,.

\noindent The structure of the semigroup $\Gamma$ was given by Patil
\& Singh (1990) under a more general assumption of almost arithmetic
sequence on the integers $m_{0} < m_{1} < \ldots < m_{p}$.
Subsequently, Patil (1993) constructed a minimal generating set
$\mathcal{G}$ for $\wp$\, which was proved to be a Gr\"{o}bner basis
by Sengupta (2003a). Al-Ayyoub(2009) pointed out a gap in Sengupta's
proof(2003a) in one particular case. However, it is still not known
whether the minimal generating set $\mathcal{G}$ for $\wp$\, in that
particular case is a Gr\"{o}bner basis with respect to some suitable
monomial order.

\noindent We restrict our attention only in the case of an
arithmetic sequence. A complete description of all the syzygy
modules was given by Sengupta (2003b), when $p=3$ and $m_{0} < m_{1}
< \ldots < m_{p}$ forms an arithmetic sequence. An explicit
description of a minimal generating set $\mathcal{G}$ for $\wp$\, in
case of arithmetic sequence is given
in Sengupta (2000)\,and in Maloo-Sengupta (2003). We recall the set putting $Y=X_{p}$ :\\[2mm]
\noindent $\mathcal{G}=\{\xi_{i\,,\,j} \mid i, j \in [1,p-2]
\}\cup\{\phi_{i} \mid i\in [0,p-2]\} \cup\{\psi_{b\,,\,j} \mid j\in
[0,p-b-1]\}\cup\{\theta\}$ where
\begin{list}{\labelitemi}{\leftmargin=1em}
\item $\xi_{i\,,\,j}=\begin{cases}X_{i}X_{j}-X_{i+j}X_{0}\hspace*{0.85in} {\rm if}
\hspace*{0.25in} i+j \leq
p-1\\[2mm]
X_{i}X_{j}-X_{i+j+1-p}X_{p-1}\hspace*{0.40in} {\rm if}
\hspace*{0.25in} i+j \geq p
\end{cases}$

\item $\phi_{i}=X_{i+1}X_{p-1}-X_{i}X_{p}$\quad \quad,\quad
\quad $\psi_{b\,,\,j}=X_{b+j}X_{p}^{a}-X_{j}X_{0}^{a+d}$\quad\quad
and \quad\quad $\theta=X_{p}^{a+1}-X_{p-b}X_{0}^{a+d}$
\end{list}

\noindent Let us construct the set $\mathcal{G}'=\{\phi(i\,,\,j)
\,\mid\, i,j \in [1,p-1]\}\cup\{\psi(b\,,\,i) \,\, \mid \,\, i \in
[0,p-b] \}$ where\\[3mm]
\noindent
$\phi(i\,,\,j)=X_{i}X_{j}-X_{\epsilon(i\,,\,j)}X_{i+j-\epsilon(i\,,\,j)}$\quad
and \quad $\psi(b\,,\,i)=X_{b+i}X_{p}^{a}-X_{i}X_{0}^{a+d}$\\[3mm]
\noindent It is easy to verify that
\begin{list}{\labelitemi}{\leftmargin=1em}
\item $\xi_{i\,,\,j}=\phi(i\,,\,j) \hspace*{1.05in} {\rm if} \hspace*{0.25in}  i,j \in [1,p-2]
\hspace*{0.30in} {\rm and} \hspace*{0.25in} i+j \leq p-1$
\item $\xi_{i\,,\,j}+\phi_{i+j-p}=\phi(i,j) \hspace*{0.45in} {\rm if} \hspace*{0.25in}  i,j \in [1,p-2]
\hspace*{0.35in} {\rm and} \hspace*{0.25in} i+j \geq p$
\item $\phi_{i}=\phi(i+1,p-1) \hspace*{0.60in} {\rm if} \hspace*{0.25in} i \in
[0,p-2]$
\item $\psi_{b\,,\,i}=\psi(b\,,\,i) \hspace*{1in} {\rm if} \hspace*{0.25in} i \in
[0,p-b-1] \hspace*{0.25in} {\rm and} \hspace*{0.25in}
\theta=\psi(b\,,\,p-b)$
\end{list}
\noindent Therefore, the set $\mathcal{G}$ is contained in the ideal
generated by the set $\mathcal{G}'$ as well as their cardinalities
are equal. Also note that $\mathcal{G}' \subseteq \wp$\,. Hence,
$\mathcal{G}'$ is a minimal generating set for the ideal $\wp$\,.

\noindent Our aim is to describe a minimal Gr\"{o}bner basis for the
symmetric algebra for $\wp$. We start by proving that $\mathcal{G}'$
is a Gr\"{o}bner basis with respect to a suitable monomial order in
Section 4. This is necessary for computing the generators of the
symmetric algebra in the subsequent sections, which is the main
theme of this article. We construct a set of linear relations among
the binomials in $\mathcal{G}'$\,. Finally we prove that the set is
a minimal Gr\"{o}bner basis for the first syzygy module of $\wp$
with respect to a suitable monomial order in Section 5. We use
numerical method to show that the leading monomials of the
generating sets indeed generate the initial ideals of the defining
ideals. The advantage of this method is that it does not require any
$S$-polynomial computation. In this regard we refer to Conca, A.,
Herzog, J., Valla, G. (1996)\,. We use the notations LT and LM to
mean leading term and leading monomial. For details of Gr\"obner
bases we refer to Eisenbud (1995) \,and Cox-Little-O'Shea (1996).
This work initiates the process of computing all relations (linear
as well as non-linear) among the binomials in $\mathcal{G}$,
culminating in a structure of the Rees algebra of $\wp$
(Mukhopadhyay \& Sengupta (2009)), which is essential for
understanding smoothness of the symbolic blow-ups of such curves
(Mukhopadhyay \& Sengupta (2009)).
\section{Gr\"{o}bner basis for \boldmath{$\wp$}}
\noindent Every monomial of\,\,R\,\, can be expressed in the form
\,\,$\displaystyle{\left(\prod_{i=1}^{p}X_{i}^{\alpha_{i}}\right)X_{0}^{\alpha_{0}}}$
\,. \,\,We identify the monomial,
$\displaystyle{\left(\prod_{i=1}^{p}X_{i}^{\alpha_{i}}\right)X_{0}^{\alpha_{0}}}$\,\,
with the ordered tuple\,\,$\displaystyle{(\alpha_{1}, \ldots,
\alpha_{p}, \alpha_{0})\in \mathbb{N}^{p+1}}$ \,.\,\, Let us define
a weight function\,\,$\omega$\\

\noindent on the monomials of\,\,R\,\,by the following\,\,:\\[2mm]
\noindent$\bullet$ \,$\omega(X_{i})=m_{i} \hspace*{1in} ;
\hspace*{0.40in} i
\in [0,p]\,.$\\

\noindent$\bullet$ \,$\omega(fg)=\omega(f)\,+\,\omega(g)
\hspace*{0.30in} ; \hspace*{0.40in} \textrm{for \, any
 \, two \, monomials} \quad f \quad \textrm{and} \quad  g \quad \textrm{of} \quad
 R\,.$\\[2mm]
\noindent We say that\quad
$\displaystyle{f=\prod_{i=1}^{p}X_{i}^{\alpha_{i}}X_{0}^{\alpha_{0}}
\,
>_{R} \,
\,\,g=\prod_{i=1}^{p}X_{i}^{\alpha'_{i}}X_{0}^{\alpha'_{0}}}$\quad
if
and only if one of the following holds\,:\\[2mm]
\noindent $\bullet$ $\omega(f) \,>\,\omega(g)$\,.\\[2mm]
\noindent $\bullet$ $\omega(f) \,=\,\omega(g)$ \, and the right-most
non-zero entry in the difference\\[2mm]
$\displaystyle{(\alpha_{1}\,-\,\alpha'_{1} \, , \, \ldots \ldots\, ,
\, \alpha_{p}\,-\,\alpha'_{p} \, , \,
\alpha_{0}\,-\,\alpha'_{0})}$\quad is negative \,.

\begin{remark}  Let\, $f$ \,and\, $g$ \,are two monomials of\, $R$\,.
One can easily check that:
\begin{list}{\labelitemi}{\leftmargin=1em}

\item $f\,-\,g \,\in \,\wp \quad \Longleftrightarrow \quad
\eta(f\,-\,g)\,=\,0 \quad \Longleftrightarrow \quad
\omega(f)=\omega(g)$

\item $\omega(X_{i}) \neq \omega(X_{j})$ for $i \neq j$ and $i,j \in [0,p]$

\item $\omega(X_{i})+\omega(X_{j})=\omega(X_{0})+\omega(X_{i+j})$ for $i+j < p$ and $i,j \in [1,p-1]$

\item $\omega(X_{i})+\omega(X_{j})=\omega(X_{p})+\omega(X_{i+j-p})$ for $i+j \leq p$ and $i,j \in [1,p-1]$
\end{list}
\end{remark}

\begin{lemma}{\it Let $m$ be the smallest integer which satisfies the
relation $mm_{p}=nm_{0}+m_{i}$ \, with \, $m\,,\,n \, \geq \, 1$
\,and \, $0 \,\leq \,i \,< \,p$ \, then : $m=a+1 $\,,\,$i=p-b$ and
$n=a+d$. }
\end{lemma}

\proof $mm_{p}=nm_{0}+m_{i} \Longrightarrow mpd=(n-m+1)m_{0}+id
\Longrightarrow mp-i=0({\rm modulo}\quad m_{o})$ since  ${\rm
gcd}(m_{0},d)=1$. There exist $q \in \mathbb{N}$ such that
$mp-i=qm_{o}$\,. Note that $q=0$ implies $mp=i$ which is absurd and
so $q \geq 1$. Therefore,  $mp-i=q(ap+b)\quad {\rm with}\quad q\geq
1$ and hence  $m \geq a+1$ since, $0\leq i <p$ and $b \in [1,p]$. At
this point note that $(a+1)m_{p}=(a+d)m_{0}+m_{p-b}$. Hence,
$m=a+1$\,,\,$n=a+d$ and $i=p-b$.\qed

\begin{corollary}{\it Let $n$ be the smallest integer which satisfies the condition
$nm_{0}=mm_{p}+m_{i}$ \, with \, $m\,,\,n \, \geq \, 1$ \,and \, $0
\,< \,i \,\leq \,p$ \, then : $n=a+d $\,,\,$i=b$ and $m=a$. }
\end{corollary}

\proof Note that $mm_{p}+m_{i}=nm_{0} \Longrightarrow
(m+1)m_{p}=(n-1)m_{0}+m_{0}+(p-i)d \Longrightarrow
(m+1)m_{p}=(n-1)m_{0}+m_{p-i}$\,. Rest of the proof follows from
Lemma 4.2.\qed

\begin{lemma}{\it If we assume that \,$m\,\neq\,n \, \geq 0$\, and \,$i\neq j \,\in [0,p] $\,
and \,$l \in \{0,p\}$\, then
\,$\omega(X_{l}^{n}X_{i})\neq\omega(X_{l}^{m}X_{j})$ \,.}
\end{lemma}

\proof $\omega(X_{l}^{n}X_{i})=\omega(X_{l}^{m}X_{j})
\Longrightarrow {\rm either} \quad sm_{l}+m_{i}=m_{j} \quad {\rm or}
\quad sm_{l}+m_{j}=m_{i}$ \,\,for some \,$s$\,. This contradicts
that $\Gamma$ is minimally generated by the set\, $S$\,. Hence the
proof\,.

\begin{lemma}{\it
If\, $g\,=\,X_{0}^{m}X_{j}$ \,with\, $ m\geq 0$  and  $j \in
[1\,,\,p-1]$ and
\begin{list}{\labelitemi}{\leftmargin=1em}
\item $f_{1}\,=\,X_{p}^{n}X_{i}$ \,with\, $i \in [1\,,\,p-1]\,\, {\rm and}
\,\, 1\leq n \leq a-1 \,$
\item $f_{2}\,=\,X_{p}^{n}X_{i}$ \,with\, $i \in [1\,,\,b-1]\,\, {\rm and}
\,\, 1\leq n \leq a \,$
\end{list}
 \, then $g-f_{q}\, \notin \, \wp$ \, for \, $q=1,2$\,.}
\end{lemma}

\proof $\omega(X_{p}^{n}X_{i})=\omega(X_{0}^{m}X_{j})$ implies that
there exist $s \geq 0$ such that either $nm_{p}=mm_{0}+sd$ \, or
$nm_{p}+sd=mm_{0}$\,. Rest of the proof follows from Lemma 4.2 and
Corollary 4.3. \qed

\begin{lemma}
{\it If\, $g\,=\,X_{0}^{m}X_{p}^{n}$ \,with\, $ m\geq 0$  \,and\, $n
\in [0\,,\,a]$  \,and\,
\begin{list}{\labelitemi}{\leftmargin=1em}
\item $f_{1}\,=\,X_{p}^{l}X_{i} \quad 1 \leq l \leq a-1 \,\, {\rm
and} \,\, i \in [1\,,\,p-1]$
\item $f_{2}\,=\,X_{0}^{l}X_{i} \quad  l \geq 0 \,\, {\rm and} \,\, i \, \in
\, [1\,,\,p-1]$
\item $f_{3}\,=\,X_{p}^{l}X_{i} \quad  1 \leq l \leq a \,\,
{\rm and} \,\, i \in [1\,,\,b-1]$
\end{list}
then $g\,-\,f_{q}\, \notin \, \wp$ \, for \,$q=1,2,3$. }
\end{lemma}

\proof The following cases will arise depending on \, $m\,\,,\,\,
n\,\,,\,\,l$ \,\,:
\begin{list}{\labelenumi}{\leftmargin=1.5em}
\item[1.] $m_{i} \in \langle \,S \setminus
m_{i} \, \rangle \,$
\item[2.] $n(ap+b) =(m+1)(ap+b)+mpd+id $
\item[3.] $m(ap+b) =(n+1)(ap+b)+id$
\end{list}

\noindent Case(i) contradicts that $\Gamma$ is minimally generated
by the set\, $S$ \,. Lemma 4.2 and Corollary 4.3 take care of
Case(ii) and Case(iii)\,. Hence the proof\,.\qed

\begin{theorem}{\it The set $\mathcal{G}'$ is a Gr\"obner
Basis for $\wp$ with respect to $>_{R}$.}
\end{theorem}

\proof If\, $f$ \,is a monomial of\, $R$ \,and\\
\noindent $f\, \, \notin \,\,{\rm LT}(\mathcal{G}')=\{X_{i}X_{j} \,
\mid  \,i\,\,j\,\,\in \,[1\,,\,p-1]\,\}\cup\{X_{b+i}X_{p}^{a} \,
\mid  \,i \in [0,p-b]\, \}$  \,then\, $f$ \, must be of the
following form\,:
\begin{list}{\labelitemi}{\leftmargin=1em}
\item $X_{p}^{m}X_{i} \,\, : \,\, 1 \leq m \leq a-1 \,\, {\rm
and} \,\, i \in [1, p-1] $

\item  $X_{0}^{m}X_{i} \,\, : \,\, i \in [1,p-1]\,\, {\rm and}
\,\, m \geq 0 \,\,{\rm is \,\, an \,\, integer} $

\item  $ X_{0}^{m}X_{p}^{n} \,\, : \,\, 0 \leq n \leq a \,\, {\rm
and} \,\, m \geq 0 \,\,{\rm is \,\, an \,\, integer} $

\item  $X_{p}^{m}X_{i} \,\, : \,\, 1 \leq m \leq a \,\, {\rm and}
\,\, i \in [1,b-1] $
\end{list}

\noindent Let \, $f\,,\,g \,\notin \,{\rm LT}(\mathcal{G}')$ are two
distinct monomials of\, $R$\,. Now apply Lemma 4.4 to Lemma 4.6 to
conclude that\, $f\,-\,g \,\notin \,\wp $. Therefore ${\rm
LT}(\wp)\,=\,\langle \,{\rm LT}(\mathcal{G}')\, \rangle$\,. Hence
the proof\,.\qed

\begin{theorem}{\it The set $\mathcal{G}'$ is a minimal Gr\"obner
Basis for $\wp$ with respect to $>_{R}$ \, .}
\end{theorem}
\noindent {\bf Proof}\quad It is enough to note that no two distinct
elements of ${\rm LT}(\mathcal{G}')$ can divide each other.\qed

\section{Symmetric Algebra}
\noindent Let\, $\widehat{R}$ \,denote the polynomial ring
\,$K[\mathbb{X}\, , \, \Psi_{b}\, , \, \Phi ]$\, with old
indeterminates \, $\mathbb{X} = \{X_{1},\ldots,X_{p},X_{0}\}$\\[2mm]
\noindent and the new ones of the set\, $\Psi_{b}=\{\Psi(b\,,\,j) \,
\mid \, j \in [0,p-b]\}$  \, and the set\, $\Phi =
\cup_{i=1}^{p-1}\Phi(\mbox{p-i})$ \,\, such that\\[2mm]
\noindent $\Phi(\mbox{i})=\{\Phi(i\,,\,i),\Phi(i-1\,,\,i),\ldots,
\Phi(1\,,\,i)\}$ \,  for
every\, $i = 1, \ldots, p-1$\,.\\[2mm]
\noindent Let \,$R[t]$\, be a polynomial ring with indeterminate $t$
\, and \,$\varphi \,: \,\widehat{R}\, \longrightarrow \,R[t]$\, is a
\, $K$-algebra homomorphism defined by
\begin{list}{\labelitemi}{\leftmargin=1em}
\item $\varphi(\Psi(b\,,\,i))=\psi(b\,,\,i) t \quad
{\rm for \,\,all} \quad i \in [0,p-b]$
\item $\varphi(\Phi(i\,,\,j))=\phi(i\,,\,j) t \quad
{\rm for \,\,all} \quad i,j \in [1,p-1]$
\item $\varphi(r)=r \quad
{\rm for \,\,all} \quad r \in R$
\end{list}
\noindent Hence the Symmetric Algebra of the ideal\, $\wp$ \, is the
polynomial ring\, $\widehat{R}/\mathcal{L}$ \,where the ideal\,
$\mathcal{L}\, \subseteq \,{\rm Kernel}(\varphi)$ \,is generated by
the set\, of polynomials which are linear with respect to the
variables\, $\Psi(b\,,\,i)$ \, and\, $\Phi(i\,,\,j)$ \,.\\[2mm]
\noindent Let us write\,\, $\mathcal{M}\,=\,\{\,\,
{\displaystyle{\sum_{
i,j}}}\,f_{i,j}\,\Phi(i\,,\,j)\,+\,{\displaystyle{\sum_{l}}}f_{l}\Psi(b\,,\,l)\quad
: \quad f_{i,j}\,\,,\,\,f_{l}\,\,\in\,\, R\}.$\\[2mm]
\noindent It is easy to check that under usual addition and scalar
multiplication $\mathcal{M} \subseteq \widehat{R}$ is a module over
$R$. The module $\mathcal{N}=\mathcal{M}\cap {\rm Kernel}(\varphi)$
is a submodule of $\mathcal{M}$ and is called the first syzygy
module of the ideal $\wp$. First agree that every monomial of
$\mathcal{M}$ is of the form $\mathbb{X}^{\alpha}\mathbf{e}$ where
either $\mathbf{e}=\Psi(b\,,\,i)$\\[2mm]
\noindent or \,$\mathbf{e}=\Phi(i\,,\,j)$. Hence, every element $H$
of $\mathcal{M}$ can be expressed uniquely as
$H={\displaystyle{\sum_{i}}}\mathbb{X}^{\alpha^{i}}\mathbf{e}_{i}$
where,
$\alpha^{i}=(\alpha_{1(i)}\,,\ldots\,,\,\alpha_{p(i)},\alpha_{0(i)})
\in \mathbb{N}^{p+1}$ and \, either \,
$\mathbf{e}_{i}=\Psi(b\,,\,l(i))$
\, or \, $\mathbf{e}_{i}=\Phi(l(i)\,,\,j(i))$ with\\[2mm]
\noindent $l(i) \,,\, j(i) \,\in \mathbb{N}$. Note that this is a
finite sum because every element $H$ of $\mathcal{M}$ is a
polynomial in $\widehat{R}$.\\[2mm]
\noindent Before we proceed further, let us record the following
remark.\\[2mm]
\noindent {\bf Remark}\quad {\it It is interesting to note that
\,$\phi(i\,,\,j)=\phi(j\,,\,i)$\, for all \,$i,j \in \mathbb{N}$.
Therefore,
$\varphi(\Phi(i\,,\,j))=\\[2mm]
\noindent \varphi(\Phi(j\,,\,i))$\, for all \,$i,j \in [1,p-1]$.
Henceforth, throughout the rest of this section we
write \,$\phi(i,j),\Phi(i,j)$\\[2mm]
\noindent to mean that $i \leq j$.}\\[2mm]
\noindent Let us define a function $\varpi$ on the set of monomials
of $\mathcal{M}$ by\\[2mm]
$\hspace*{0.80in}\varpi(\mathbb{X}^{\alpha}\Psi(b\,,\,i)) \,=\,
\mathbb{X}^{\alpha}X_{p}^{a}X_{b+i}  \hspace*{0.30in} ;
\hspace*{0.30in} {\rm with} \hspace*{0.50in}    i \,\in \,[0\,,\,
p-b]$\\[2mm]
\noindent and $\hspace*{0.50in}
\varpi(\mathbb{X}^{\alpha}\Phi(i\,,\,j)) \,=\,
\mathbb{X}^{\alpha}X_{i}X_{j} \hspace*{0.50in}; \hspace*{0.30in}{\rm
with}   \hspace*{0.50in} i \,, \,j \, \in \,
[1\,,\,p-1] \hspace*{0.30in}.$\\[2mm]
\noindent We now define a monomial order\, $>_{\mathcal{M}}$\,\,on
\,$\mathcal{M}$\,\,by \,$\mathbb{X}^{\alpha}\mathbf{e} \,
>_{\mathcal{M}} \,\mathbb{X}^{\alpha'}\mathbf{e}'$ \,
iff one of the following holds \, :
\begin{list}{\labelitemi}{\leftmargin=1em}
\item[1.] $\varpi(\mathbb{X}^{\alpha}\mathbf{e})
\quad  >_{R} \quad \varpi(\mathbb{X}^{\alpha'}\mathbf{e}')$

\item[2.] $\varpi(\mathbb{X}^{\alpha}\mathbf{e})
\,=\,\varpi(\mathbb{X}^{\alpha'}\mathbf{e}')$\,, and  one of the
following holds \, :

\begin{list}{\labelitemi}{\leftmargin=1em}
\item $\mathbf{e}\,=\,\mathbb{X}^{\alpha}\Psi(b\,,\,i)$ \,and\,
$\mathbf{e}'\,=\,\mathbb{X}^{\alpha'}\Psi(b\,,\,j)$
and $i\,<\,j$

\item $\mathbf{e}=\mathbb{X}^{\alpha}\Psi(b\,,\,i)$ \,and\,
$\mathbf{e}'\,=\,\mathbb{X}^{\alpha}\Phi(l\,,\,j)$

\item $\mathbf{e}\,=\,\mathbb{X}^{\alpha}\Phi(i\,,\,j)$ \,and\,
$\mathbf{e}'\,=\,\mathbb{X}^{\alpha}\Phi(i'\,,\,j')$
\,and \, either \, $j\,>\,j'$ \, or \, $j\,=\,j'$ \,and\, $i\,>\,i'$
\end{list}
\end{list}
\begin{remark}
It is interesting to note that $\varpi(F)={\rm LM}(\varphi(F))$ for
every monomial $F$ of $\mathcal{M}$\,.
\end{remark}
\begin{definition}{\rm
An element \,$H\neq 0$ \,in \,$\mathcal{M}$ is called a {\it
relation} if \,$\varphi(H) = 0$. A {\it relation} $H$ is called a
{\it reduced\,\, relation} if \,$H=FG$ \,with \,$F \in R$ \,,\, $G
\in \mathcal{M}$ and \,$\varphi(G)=0$ implies that $F =1$. A {\it
reduced\,\,relation} $H$ is called a {\it basic\,\, relation} if
\,$H=F+G$ \,with \,$\varphi(F)=0$ \,and \,$\varphi(G)=0$ \,implies
that either \, $F=0$ \, or \, $G=0$.}
\end{definition}

\begin{remark}
Note that, by the definition of $\varphi$\,, every relation can be
written as a finite linear combination of {\it basic\,\,relations}
over $R$. Therefore, we may assume that the generators of the module
$\mathcal{N}$ over $R$ are {\it basic\,\,relations}. Henceforth, the
term {\it relation} stands for a {\it basic\,\,relation}.
\end{remark}

\noindent For our convenience from now on we will use underlined
terms to represent leading terms. Unless otherwise specified, the
symbols $H_{l} , c_{l} $ for each integer $l$ will represent a
monomial of $\mathcal{M}$ and an element of $K$ respectively.
Throughout the rest of this article, $\alpha_{0} \,,\, \alpha_{p} \,
\in \, \mathbb{N}$ and $\alpha\,,\,\beta \, \in \,
\mathbb{N}^{p+1}$\,.
\begin{lemma} If $H=\underline{c_{1}H_{1}}+{\displaystyle{\sum_{l}}} c_{l}H_{l}$
is a relation, then there exist $l'$ such that ${\rm
LM}(\varphi(H_{1}))={\rm LM}(\varphi(H_{l'}))$ with
$H_{l'}\,\neq\,H_{1}$\,.
\end{lemma}

\proof Follows from definition of $\varphi$ and monomial order
$>_{\mathcal{M}}$\,.\qed
\begin{lemma} If $H=\underline{c_{1}\mathbb{X}^{\alpha}\Psi(b\,,\,i)}+{\displaystyle{\sum_{l}}}
c_{l}H_{l}$ is a relation then
$\mathbb{X}^{\alpha}\,\neq\,X_{0}^{\alpha_{0}}$\,.
\end{lemma}

\proof  If possible assume that
$H=\underline{c_{1}X_{0}^{\alpha_{0}}\Psi(b\,,\,i)}+{\displaystyle{\sum_{l}}}
c_{l}H_{l}$. According to Lemma 5.4.,
 there exist\\[2mm] \noindent $l'$
such that \, ${\rm LM}(\varphi(H_{l'}))\,=\,{\rm
LM}(\varphi(X_{0}^{\alpha_{0}}\Psi(b\,,\,i)))\,=\,X_{0}^{\alpha_{0}}X_{p}^{a}X_{b+i}t$.
It is clear from the explicit \\[2mm]
\noindent description of ${\rm LT}(\mathcal{G}')$, that no such
$H_{l'}$ exist. Hence the proof\,.\qed

\begin{lemma} If
$H=\underline{c_{1}\mathbb{X}^{\alpha}\Phi(i\,,\,j)}+{\displaystyle{\sum_{l}}}
c_{l}H_{l}$ is a relation then
$\mathbb{X}^{\alpha}\,\neq\,X_{0}^{\alpha_{0}}$\,.
\end{lemma}
\proof Similar to Lemma 5.5.\qed

\begin{lemma}{\it If
$H=\underline{c_{1}\mathbb{X}^{\alpha}\Psi(b\,,\,p-b)}+{\displaystyle{\sum_{l}}}
c_{l}H_{l}$ is a relation then\\

\noindent $\mathbb{X}^{\alpha}\,\notin \{X_{i}^{\alpha_{i}}X_{j}^{q}
\,:\,j\in [0,p] \,\,{\rm and} \,\, i \in \{0\,,\,p\} \,\,{\rm and}
\,\, q \in \{0,1\} \,\}$.}
\end{lemma}

\proof If possible assume that
$H=\underline{c_{1}X_{i}^{\alpha_{i}}X_{j}^{q}\Psi(b\,,\,p-b)}+{\displaystyle{\sum_{l}}}
c_{l}H_{l}$. According to Lemma 5.4.,\\

\noindent there exist $l'$ such that ${\rm LM}(\varphi(H_{l'}))={\rm
LM}(\varphi(X_{i}^{\alpha_{i}}X_{j}^{q}\Psi(b\,,\,p-b)))=X_{i}^{\alpha_{i}}X_{j}^{q}X_{p}^{a+1}t$.
From the explicit\\

\noindent description of ${\rm LT}(\mathcal{G}')$, it is clear that
for $q=0$ or, for $j \notin [b,p-1]$ no such $H_{l'}$ exist and
for\\

\noindent $j \in [b,p-1]$ one can write
$H_{l'}=X_{i}^{\alpha_{i}}X_{p}\Psi(b\,,\,j-b)$ but then ${\rm
LM}(H)=H_{l'}$ will contradict the leading\\

\noindent monomial assumption. Hence the proof.\qed

\begin{lemma}{\it If
$H=\underline{c_{1}\mathbb{X}^{\alpha}\Phi(i\,,\,j)}+{\displaystyle{\sum_{l}}}
c_{l}H_{l}$ is a {\it relation} then $X_{l'} \mid
\mathbb{X}^{\alpha}$ for some $0<l' <j$.}
\end{lemma}

\proof If possible assume that
$H=\underline{c_{1}\mathbb{X}^{\alpha}\Phi(i\,,\,j)}+{\displaystyle{\sum_{l}}}
c_{l}H_{l}$ and $X_{l'} \nmid \mathbb{X}^{\alpha}$ where $0<l'
<j$.\\[1mm]

\noindent According to Lemma 5.4., there exist $q$ such that ${\rm
LM}(\varphi(H_{q}))={\rm
LM}(\varphi(\mathbb{X}^{\alpha}\Phi(i\,,\,j)))=
\mathbb{X}^{\alpha}X_{i}X_{j}t$.\\[1mm]

\noindent From the explicit description of ${\rm LT}(\mathcal{G}')$,
it is clear that $H_{q}=\begin{cases}
\mathbb{X}^{\beta}\Phi(i\,,\,l) \quad
{\rm with} \quad l>j\\[2mm]
\mathbb{X}^{\beta}\Phi(j\,,\,l) \quad
{\rm with} \quad l\geq j\\[2mm]
\mathbb{X}^{\beta}\Phi(s\,,\,l) \quad
{\rm with} \quad s,l \,\geq j\\[2mm]
\end{cases}$\\[1mm]

\noindent Therefore in all the cases ${\rm LM}(H)=H_{q}$. This
contradiction proves the result.\qed

\section{Gr\"obner basis for the first syzygy module}
\noindent For systematic reason set \quad
$\Phi(i\,,\,j)=\Phi(j\,,\,i)$ \quad {\rm and}  \quad
$\Phi(i\,,\,j)=0 \quad {\rm if} \quad i,j \notin
[1,p-1]$\,.\\
\noindent Let us construct the set $\widehat{\mathcal{G}}$ whose
elements are the following(with underlined leading terms)\,:
\begin{list}{\labelitemi}{\leftmargin=1em}
\item $A(i\,;\,b\,,\,j)=
\underline{X_{i}\Psi(b\,,\,j)}-X_{b+i+j-\epsilon(i\,,\,
b+j)}\Psi(b,\epsilon(i\,,\,b+j)-b)
-X_{p}^{a}\Phi(i\,,\,b+j)\\[2mm]
\hspace*{0.70in}+X_{0}^{a+d}[\Phi(i\,,\,j)-\Phi(b+i+j-p\,,\,p-b)] $

\item $B(i\,,\,j)=
\underline{X_{i}X_{j}\Psi(b\,,\,p-b)}-X_{\epsilon(i\,,\,j)}X_{i+j-\epsilon(i\,,\,j)}\Psi(b\,,\,p-b)
-\psi(b\,,\,p-b)\Phi(i\,,\,j)$

\item $L(l\,;\,i\,,\,j)=
\underline{X_{l}\Phi(i\,,\,j)}\, - \, X_{j}\Phi(i\,,\,l)
+X_{\tau(i\,,\,j)}\Phi(i+j-\tau(i\,,\,j)\,,\,l)
-X_{\tau(i\,,\,l)}\Phi(i+l-\tau(i\,,\,l)\,,\,j)$
\end{list}

\noindent Our aim is to prove that the set\\[2mm]
$\widehat{\mathcal{G}}=\{A(i\,;\,b\,,\,j)\,\mid\, i\in[1,p]\quad
{\rm and}\quad j\in [0,p-b-1]\}\,\cup\,\{B(i\,,\,j)\,\,\mid\,\, i,j
\in
[1,p-1]\}\\[2mm]
\hspace*{0.15in}\cup\{L(l\,;\,i\,,\,j)\,\, \mid\,\, l,i,j \in
[1,p-1]\quad {\rm
with}\quad i\leq j \quad {\rm and} \quad l < j \}$\\[2mm]
\noindent is a minimal Gr\"obner basis of \,$\mathcal{N}$. \,Note
that \,$\widehat{\mathcal{G}}\,\subseteq \, \mathcal{N}$.
\begin{theorem}{\it
$\widehat{\mathcal{G}}$\, is a Gr\"obner basis for the first Syzygy
module of $\wp$ with respect to $\displaystyle{>_{\mathcal{M}}}$.}
\end{theorem}

\proof If $f$ is a monomial of $\mathcal{M}$ and $f \notin {\rm
LT}(\widehat{\mathcal{G}})$ where,\\[2mm]
${\rm
LT}(\widehat{\mathcal{G}})=\{X_{i}X_{j}\Psi(b\,,\,p-b)\,\,\mid\,\,
i,j \in
[1,p-1]\}\,\cup\,\{X_{l}\Phi(i\,,\,j)\,\,\mid\,\,  l<j\,\,{\rm and}\,\,i \leq j\}\\[2mm]
\hspace*{0.40in}\,\,\cup\,\{X_{l}\Psi(b\,,\,j)\,\,\mid\,\,
 l\,\in\,[1,p]\,\,\,
{\rm and}\,\,\,j\,\in \,[0,p-b-1]\}$\\[2mm]
then \,$f$ \,must be one of the following \,:\\[2mm]
\noindent Case(i)\,:\,
$X_{0}^{\alpha_{0}}\Psi(b\,,\,i)$\hspace*{0.50in} Case(ii)\,:\,
$X_{0}^{\alpha_{0}}X_{i}\Psi(b\,,\,p-b)$\hspace*{0.50in}
Case(iii)\,:\,
$X_{p}^{\alpha_{p}}X_{i}\Psi(b\,,\,p-b)$\\[2mm]
\noindent Case(iv)\,:\, $\mathbb{X}^{\alpha}\Phi(i\,,\,j)\quad {\rm
with} \quad X_{l}\,\nmid \, \mathbb{X}^{\alpha} \quad {\rm for}
\quad
0\,<\,l<\,j$\\[2mm]
Rest of the proof follows from Lemma 5.5. to Lemma 5.8.\qed

\begin{theorem}{\it The set
$\widehat{\mathcal{G}}$ is a minimal Gr\"{o}bner basis for the first
Syzygy module of \, $\wp$.}
\end{theorem}

\proof It is enough to note that no two distinct elements of ${\rm
LT}(\widehat{\mathcal{G}})$ can divide each other.\qed

\bigskip

\Acknow  I would like to express my heartfelt gratitude and sincere
thanks to my supervisor Prof. Indranath Sengupta for his guidance,
help and encouragement. I am also grateful and thankful to Prof.
Tarun Kumar Mukherjee  and Prof. Himadri  Sarkar for their help and
encouragement.

\begin{flushleft}
\Large \bf References
\end{flushleft}

{\small \noindent[1]\,\, Conca, A., Herzog, J., Valla, G. (1996) .
Sagbi bases with application to blow-up algebras.
\emph{J.reine.angew.Math} 474: 113-138.

\noindent[2]\,\,Eisenbud, E. (1995). \emph{Commutative Algebra with
a View Toward Algebraic Geometry}. New York: Springer - Verlag.

\noindent[3]\,\, Cox, D., Little, J., O'Shea, D. (1996).
\emph{Ideals, Varieties and Algorithms}. New York: Springer-Verlag.

\noindent[4]\,\, Patil, D. P., Singh, B. (1990). Generators for the
derivation modules and the relation ideals of certain curves.
\emph{Manuscripta Math.} 68: 327-335.

\noindent[5]\,\, Patil, D. P. (1993). Minimal sets of generators for
the relation ideal of certain monomial curves. \emph{Manuscripta
Math.} 80: 239-248.

\noindent[6]\,\, Sengupta, I. (2003a). A Gr\"obner bases for certain
affine monomial curves. \emph{Communications in Algebra} 31(3):
1113-1129.

\noindent[7]\,\, Sengupta, I. (2003b). A minimal free resolution for
certain monomial curves in $\mathbb{A}^{4}$. \emph{Communications in
Algebra} 31(6): 2791-2809.

\noindent[8]\,\, Sengupta, I. (2000). \emph{Betti Numbers,
Gr\"{o}bner  Basis and Syzygies for Certain Affine Monomial Curves.}
Thesis, Indian Institute of Science, Bangalore, India.

\noindent[9]\,\, Maloo, A.K., Sengupta, I. (2003). Criterion for
Complete Intersection for Certain Monomial Curves. \emph{Advances in
Algebra and Geometry}, University of Hyderabad Conference 2001,
Edited by C.Musili, Hindustan Book Agency, pp. 179-184.

\noindent[10]\,\,Al-Ayyoub,I.(2009). Reduced Gr\"obner bases of
certain toric varieties; A new short proof. \emph{Communications in
Algebra} 37(9): 2945-2955.

\noindent[11]\,\, Mukhopadhyay, D., Sengupta, I. (2009). The Rees
Algebra for Certain Monomial Curves. \emph{Preprint.}

\noindent[12]\,\, Mukhopadhyay, D., Sengupta, I. (2009). On The
Smoothness of Blowups for Certain Monomial Curves. \emph{Preprint.}
}

\medskip

\noindent Address : Acharya Girish Chandra Bose College, 35, Scott
Lane, Kolkata, WB 700009, INDIA.
\medskip
\noindent E-mail : mdebasish01@yahoo.co.in

\end{document}